\documentclass{amsart}
\usepackage{amsmath,amssymb}
\usepackage[dvips]{graphics}
\usepackage[all]{xy}
\newtheorem{Theorem}{Theorem}[section]
\newtheorem{Lemma}[Theorem]{Lemma}
\newtheorem{Proposition}[Theorem]{Proposition}

\theoremstyle{Definition}

\newtheorem{Example}[Theorem]{Example}

\theoremstyle{Remark}
\newtheorem{Remark}[Theorem]{Remark}

\makeatletter
\@addtoreset{figure}{section}%{subsection}
\def\@thmcountersep{-}
\makeatother

\numberwithin{equation}{section}

\begin{document} 

\title{Homotopy on spatial graphs and generalized Sato-Levine invariants}

%    Information for first author
\author{Ryo Nikkuni}
\address{Department of Mathematics, Faculty of Education, Kanazawa University, Kakuma-machi, Kanazawa, Ishikawa, 920-1192, Japan}
\email{nick@ed.kanazawa-u.ac.jp}
%\thanks{}

%    General info

\date{}

\dedicatory{}


\begin{abstract}
Edge-homotopy and vertex-homotopy are equivalence relations on spatial graphs which are generalizations of Milnor's link-homotopy. Fleming and the author introduced some edge (resp. vertex)-homotopy invariants of spatial graphs by applying the Sato-Levine invariant for the constituent 2-component algebraically split links. In this paper, we construct some new edge (resp. vertex)-homotopy invariants of spatial graphs without any restriction of linking numbers of the constituent 2-component links by applying the generalized Sato-Levine invariant. 
\end{abstract}

\maketitle

\section{Introduction} 

Throughout this paper we work in the piecewise linear category. Let $G$ be a finite graph. An embedding $f$ of $G$ into the $3$-sphere ${\mathbb S}^{3}$ is called  a {\it spatial embedding of $G$} or simply a {\it spatial graph}. We call the image of $f$ restricted on a cycle (resp. mutually disjoint cycles) in $G$ a {\it constituent knot} (resp. {\it constituent link}) of $f$, where a {\it cycle} is a graph homeomorphic to a circle. A spatial embedding of a planar graph is said to be {\it trivial} if it is ambient isotopic to an embedding of the graph into a $2$-sphere in ${\mathbb S}^{3}$. A spatial embedding $f$ of $G$ is said to be {\it split} if there exists a $2$-sphere $S$ in ${\mathbb S}^{3}$ such that $S\cap f(G)=\emptyset$ and each connected component of ${\mathbb S}^{3}-S$ has intersection with $f(G)$, and otherwise $f$ is said to be {\it non-splittable}. 

\begin{figure}[htbp]
      \begin{center}
\scalebox{0.34}{\includegraphics*{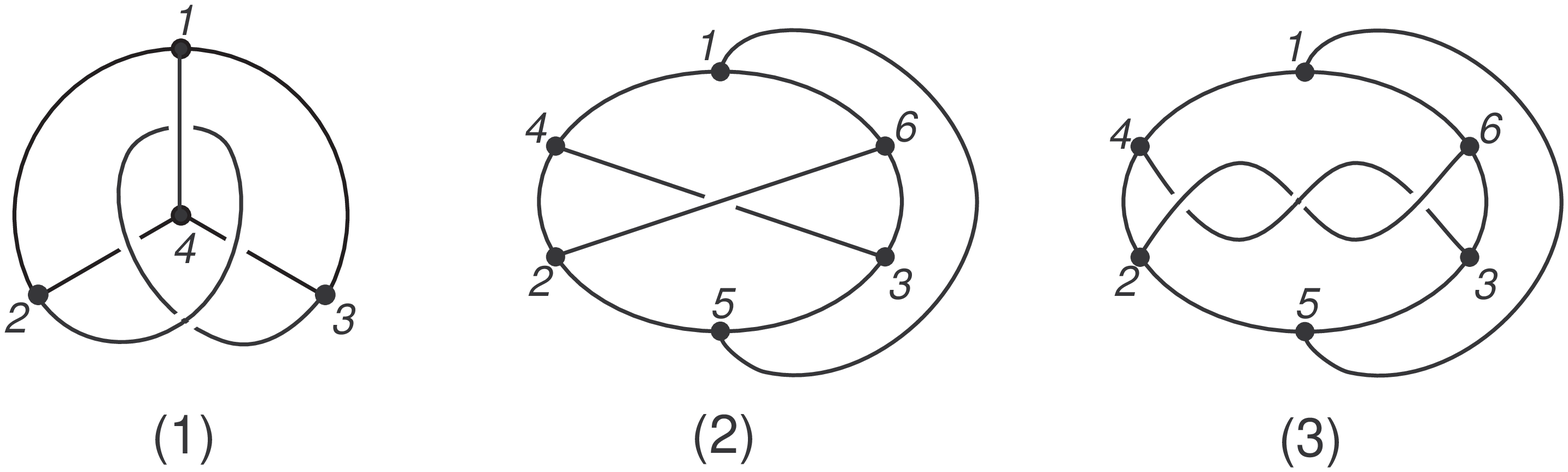}}
      \end{center}
   \caption{}
  \label{knownEx1}
\end{figure}

\begin{figure}[htbp]
      \begin{center}
\scalebox{0.34}{\includegraphics*{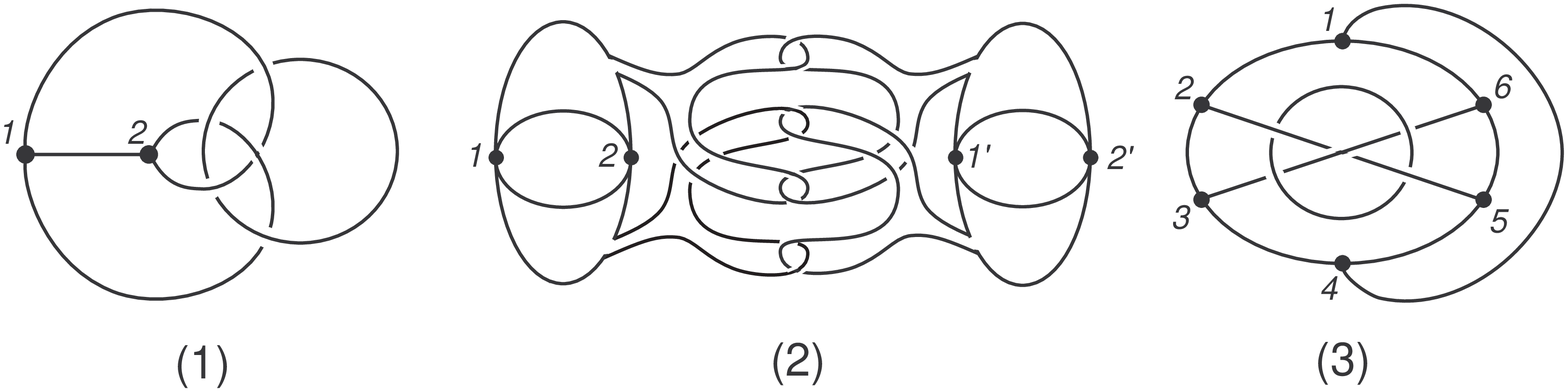}}
      \end{center}
   \caption{}
  \label{knownEx2}
\end{figure}

Two spatial embeddings of $G$ are said to be {\it edge-homotopic} if they are transformed into each other by {\it self crossing changes} and ambient isotopies, where a self crossing change is a crossing change on the same spatial edge, and {\it vertex-homotopic} if they are transformed into each other by crossing changes on two adjacent spatial edges and ambient isotopies. These equivalence relations were introduced by Taniyama \cite{taniyama94b} as generalizations of Milnor's {\it link-homotopy} on oriented links \cite{milnor54}, namely if $G$ is a mutually disjoint union of cycles then these are none other than link-homotopy. It is known that edge (resp. vertex)-homotopy on spatial graphs behaves quite differently than link-homotopy on oriented links. Taniyama introduced the {\it $\alpha$-invariant} of spatial graphs by taking a weighted sum of the second coefficient of the Conway polynomial of the constituent knots \cite{taniyama94a}. By applying the $\alpha$-invariant, it is shown that the spatial embedding of $K_{4}$ as illustrated in Fig. \ref{knownEx1} (1) is not trivial up to edge-homotopy, and two spatial embeddings of $K_{3,3}$ as illustrated in Fig. \ref{knownEx1} (2) and (3) are not vertex-homotopic. Note that each of these spatial graphs does not have a consituent link. On the other hand, some invariants of spatial graphs defined by taking a weighted sum of the third coefficient of the Conway polynomial of the constituent $2$-component links were introduced by Taniyama as ${\mathbb Z}_{2}$-valued invariants if the linking numbers are even \cite{taniyama0}, and by Fleming and the author as integer-valued invariants if the linking numbers vanish \cite{fleming-nikkuni05}. By applying these invariants, it is shown that each of the spatial graphs as illustrated in Fig. \ref{knownEx2} (1) and (2) is non-splittable up to edge-homotopy, and the spatial graph as illustrated in Fig. \ref{knownEx2} (3) is non-splittable up to vertex-homotopy. Note that each of these spatial graphs does not contain a constituent link which is not trivial up to link-homotopy. 

Our purpose in this paper is to construct some new edge (resp. vertex)-homotopy invariants of spatial graphs without any restriction of linking numbers of the constituent 2-component links by applying a weighted sum of the {\it generalized Sato-Levine invariant}. Here the generalized Sato-Levine invariant $\tilde{\beta}(L)=\tilde{\beta}(K_{1},K_{2})$ is an ambient isotopy invariant of an oriented $2$-component link $L=K_{1}\cup K_{2}$ which appears in various ways independently \cite{akhmetiev98}, \cite{akhmetiev-repovs98}, \cite{kirk-livingston97}, \cite{livingston03}, \cite{kanenobu-miyazawa-tani98}, \cite{nakanishi02} and can be calculated by  
\begin{eqnarray*}
\tilde{\beta}(L)=a_{3}(L)-{\rm lk}(L)\left\{a_{2}(K_{1})+a_{2}(K_{2})\right\},
\end{eqnarray*}
where $a_{i}$ denotes the $i$-th coefficient of the Conway polynomial and ${\rm lk}(L)={\rm lk}(K_{1},K_{2})$ denotes the linking number of $L$. It is known that the original {\it Sato-Levine invariant} $\beta(L)$ \cite{sato84} coincides with $a_{3}(L)$ if ${\rm lk}(L)=0$ \cite{cochran85}, \cite{sturm90}. Thus in this case we have that $\tilde{\beta}(L)={\beta}(L)$. As a consequence, our invariants are generalizations of Fleming and the author's homotopy invariants of spatial graphs defined in \cite{fleming-nikkuni05}. 

This paper is organized as follows. In the next section, we show some formulas about the generalized Sato-Levine invariant of oriented $2$-component links needed later. In section $3$, we give the definitions of our invariants and state their invariance up to edge (resp. vertex)-homotopy. In section $4$, we give some examples. 

\section{Some formulas about the generalized Sato-Levine invariant} 

We first show the change in the generalized Sato-Levine invariant of oriented $2$-component links which differ by a single self crossing change. 

\begin{Lemma}\label{keylemma}
Let $L_{+}=J_{+}\cup K$ and $L_{-}=J_{-}\cup K$ be two oriented $2$-component links and $L_{0}=J_{1}\cup J_{2}\cup K$ an oriented $3$-component link which are identical except inside the depicted regions as illustrated in Fig. \ref{skein_sc}. Suppose that ${\rm lk}(L_{+})={\rm lk}(L_{-})=m$. Then it holds that 
\begin{eqnarray*}
\tilde{\beta}(L_{+})-\tilde{\beta}(L_{-})
={\rm lk}(K, J_{i})\left\{m-{\rm lk}(K, J_{i})\right\}\ (i=1,2). 
\end{eqnarray*}
\begin{figure}[htbp]
      \begin{center}
\scalebox{0.4}{\includegraphics*{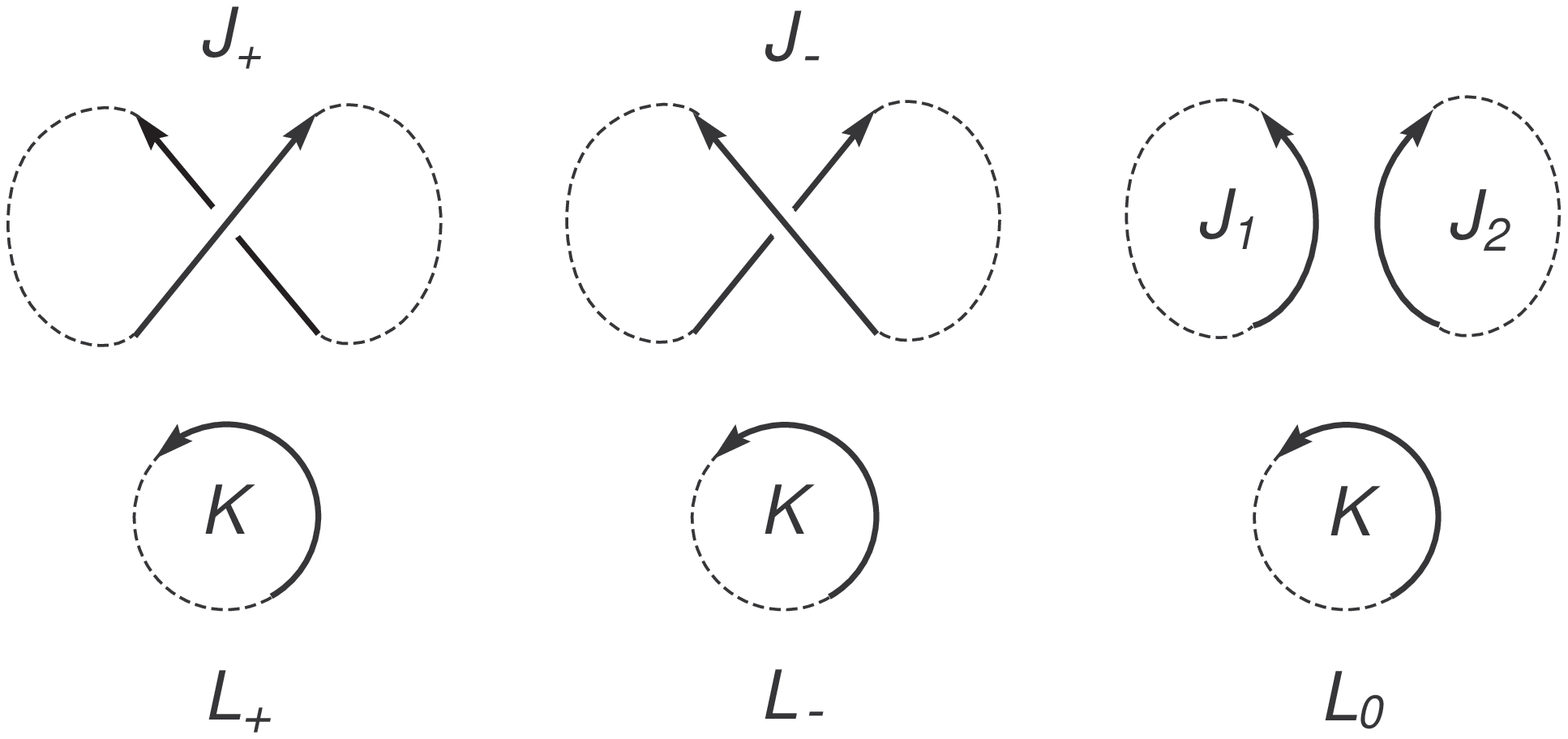}}
      \end{center}
   \caption{}
  \label{skein_sc}
\end{figure} 
\end{Lemma}

We remark here that this formula has already known (see \cite[Theorem 8.7]{livingston03} for example), but we give the proof again for readers' convenience. 

\begin{proof} 
By the skein relation of the Conway polynomial and well-known formulas for the first coefficient of the Conway polynomial of an oriented $2$-component link (cf. \cite{kauffman83}) and for the second coefficient of the Conway polynomial of an oriented $3$-component link (cf. \cite{hosokawa58}, \cite{hartley83}, \cite{hoste85}), we have that 
\begin{eqnarray}
\label{skein0} a_{2}(J_{+})-a_{2}(J_{-})&=&{\rm lk}(J_{1},J_{2}), \\
\label{skein1} a_{3}(L_{+})-a_{3}(L_{-})
&=&{\rm lk}(J_{1},J_{2}){\rm lk}(J_{2},K)+{\rm lk}(J_{2},K){\rm lk}(J_{1},K)\\
&&+{\rm lk}(J_{1},K){\rm lk}(J_{1},J_{2}). \nonumber
\end{eqnarray}
Note that 
\begin{eqnarray}\label{lk}
{\rm lk}(J_{1},K)+{\rm lk}(J_{2},K)=m. 
\end{eqnarray}
Thus by (\ref{skein0}), (\ref{skein1}) and (\ref{lk}), we have that 
\begin{eqnarray*}
\tilde{\beta}(L_{+})-\tilde{\beta}(L_{-})&=&
a_{3}(L_{+})-m\left\{a_{2}(J_{+})+a_{2}(K)\right\}
-a_{3}(L_{-})+m\left\{a_{2}(J_{-})+a_{2}(K)\right\}\\
&=&a_{3}(L_{+})-a_{3}(L_{-})-m\left\{a_{2}(J_{+})-a_{2}(J_{-})\right\}\\
&=&{\rm lk}(J_{1},J_{2})\left\{m-{\rm lk}(J_{1},K)\right\}+{\rm lk}(J_{2},K){\rm lk}(J_{1},K)\\
&&+{\rm lk}(J_{1},K){\rm lk}(J_{1},J_{2})-m{\rm lk}(J_{1},J_{2})\\
&=&{\rm lk}(J_{2},K){\rm lk}(J_{1},K). 
\end{eqnarray*}
Therefore by (\ref{lk}) we have the result. 
\end{proof}

Next we investigate the change in the generalized Sato-Levine invariant of oriented $2$-component links which differ by inverting the orientation on one of the components. The original definition implies that the value of the Sato-Levine invariant does not depend on the orientation of each component. But the value of the generalized Sato-Levine invariant depends on it in general. 

\begin{Theorem}\label{gsl_ori}
Let $L=J_{1}\cup J_{2}$ be an oriented $2$-component link with ${\rm lk}(L)=m$. Let $L'=(-J_{1})\cup J_{2}$ be the oriented $2$-component link obtained from $L$ by inverting the orientation of $J_{1}$. Then it holds that 
\begin{eqnarray*}
\tilde{\beta}(L)-\tilde{\beta}(L')
=\frac{1}{6}(m^{3}-m). 
\end{eqnarray*} 
\end{Theorem}

\begin{proof}
Let $T_{m}$ and $T'_{m}$ be two oriented $2$-component links as illustrated in Fig. \ref{(2,m)torus}. By a direct calculation we have that ${\rm lk}(T_{m})=m$, ${\rm lk}(T'_{m})=-m$ and 
\begin{eqnarray}
\label{beta_cal1} a_{3}(T_{m})&=&\frac{1}{6}(m^{3}-m),\\ 
\label{beta_cal2} a_{3}(T'_{m})&=&0.
\end{eqnarray} 
\begin{figure}[htbp]
      \begin{center}
\scalebox{0.5}{\includegraphics*{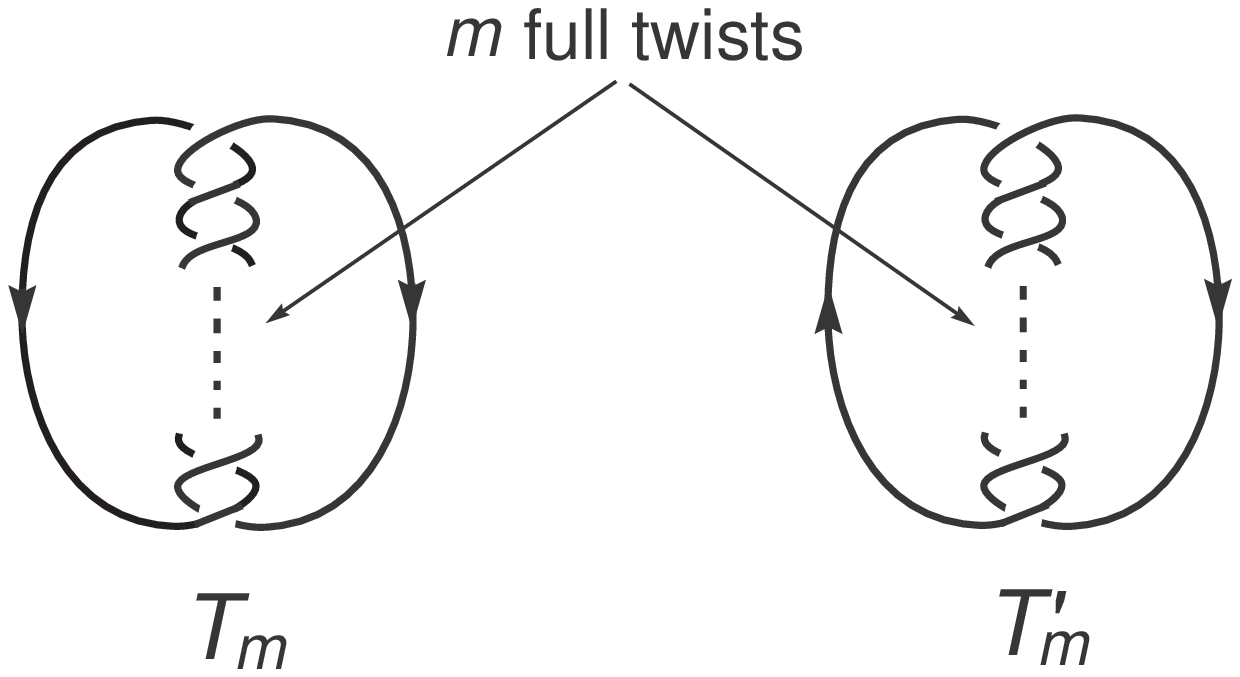}}
      \end{center}
   \caption{}
  \label{(2,m)torus}
\end{figure} 

By the classification of oriented $2$-component links up to link-homotopy \cite{milnor54}, we have that $L$ and $T_{m}$ are transformed into each other by self crossing changes and ambient isotopies. Namely there exists a sequence of oriented $2$-component links $L=L_{0},L_{1},\ldots,L_{k-1},L_{k}=T_{m}$ such that $L_{i}$ is obtained from $L_{i-1}$ by a single self crossing change $c_{i}$ ($i=1,2,\ldots,k$). Let $M_{i}=K_{0}^{(i)}\cup K_{1}^{(i)}\cup K_{2}^{(i)}$ be the oriented $3$-component link obtained from $L_{i-1}$ by smoothing it on $c_{i}$, where $K_{1}^{(i)}$ and $K_{2}^{(i)}$ are actual separated parts ($i=1,2,\ldots,k$). Namely we have the upper part of a skein tree as follows. 
\begin{eqnarray}\label{skeintree1}
\xymatrix{
L \ar[dr] \ar[r]^{c_{1}} & L_{1} \ar[dr] \ar[r]^{c_{2}} & L_{2} \ar[dr] \ar[r]^{c_{3}} & \cdots \ar[dr] \ar[r]^{c_{k-1}} & L_{k-1} \ar[dr] \ar[r]^{c_{k}}& T_{m}\\
& M_{1} & M_{2} & \cdots & M_{k-1} & M_{k} &\\
}
\end{eqnarray}
We define a sign of $c_{i}$ as follows: $\varepsilon(c_{i})=1$ if $c_{i}$ changes a positive crossing into a negative crossing and $-1$ if $c_{i}$ changes a negative crossing into a positive crossing. Since ${\rm lk}(L_{i})=m\ (i=1,\ldots,k)$, we have that 
\begin{eqnarray}\label{Mlk}
{\rm lk}(K_{0}^{(i)},K_{1}^{(i)})+{\rm lk}(K_{0}^{(i)},K_{2}^{(i)})=m 
\end{eqnarray}
for $i=1,2,\ldots,k$. Then by (\ref{beta_cal1}), (\ref{skeintree1}) and (\ref{Mlk}), we have that 
\begin{eqnarray}
a_{3}(L)&=&a_{3}(T_{m})+\sum_{i=1}^{k}\varepsilon(c_{i})a_{2}(M_{i}) \nonumber\\
&=&\frac{1}{6}(m^{3}-m)
+\sum_{i=1}^{k}\varepsilon(c_{i})
\Big\{
{\rm lk}(K_{0}^{(i)},K_{1}^{(i)}){\rm lk}(K_{1}^{(i)},K_{2}^{(i)})\nonumber\\
&&+{\rm lk}(K_{1}^{(i)},K_{2}^{(i)}){\rm lk}(K_{0}^{(i)},K_{2}^{(i)})
+{\rm lk}(K_{0}^{(i)},K_{2}^{(i)}){\rm lk}(K_{0}^{(i)},K_{1}^{(i)})
\Big\}\nonumber\\
&=&\frac{1}{6}(m^{3}-m)+\sum_{i=1}^{k}\varepsilon(c_{i})
\Big(
\left\{m-{\rm lk}(K_{0}^{(i)},K_{2}^{(i)})\right\}{\rm lk}(K_{1}^{(i)},K_{2}^{(i)})\nonumber\\
&&+{\rm lk}(K_{1}^{(i)},K_{2}^{(i)}){\rm lk}(K_{0}^{(i)},K_{2}^{(i)})
+{\rm lk}(K_{0}^{(i)},K_{2}^{(i)}){\rm lk}(K_{0}^{(i)},K_{1}^{(i)})
\Big)\nonumber\\
&=&
\frac{1}{6}(m^{3}-m)+\sum_{i=1}^{k}\varepsilon(c_{i})
\Big\{m{\rm lk}(K_{1}^{(i)},K_{2}^{(i)})
+{\rm lk}(K_{0}^{(i)},K_{2}^{(i)}){\rm lk}(K_{0}^{(i)},K_{1}^{(i)})
\Big\}. \label{a3_1}
\end{eqnarray}

In the same way as above, we have that $L'$ and $T'_{m}$ are transformed into each other by self crossing changes and ambient isotopies and obtain the following skein tree from (\ref{skeintree1}) immediately: 
\begin{eqnarray}\label{skeintree2}
\xymatrix{
L' \ar[dr] \ar[r]^{c'_{1}} & L'_{1} \ar[dr] \ar[r]^{c'_{2}} & L'_{2} \ar[dr] \ar[r]^{c'_{3}} & \cdots \ar[dr] \ar[r]^{c'_{k-1}} & L'_{k-1} \ar[dr] \ar[r]^{c'_{k}}& T'_{m}\\
& M'_{1} & M'_{2} & \cdots & M'_{k-1} & M'_{k} &\\
}
\end{eqnarray}
We also denote $M'_{i}={K'}_{0}^{(i)}\cup {K'}_{1}^{(i)}\cup {K'}_{2}^{(i)}$, where ${K'}_{1}^{(i)}$ and ${K'}_{2}^{(i)}$ are actual separated parts obtained by smoothing $L'_{i-1}$ on $c'_{i}$ ($i=1,2,\ldots,k$). Note that 
\begin{eqnarray}
\label{key0} \varepsilon(c_{i})&=&\varepsilon(c'_{i}), \\
\label{key1} {\rm lk}(K_{1}^{(i)},K_{2}^{(i)})&=&{\rm lk}({K'}_{1}^{(i)},{K'}_{2}^{(i)}),\\
\label{key2} {\rm lk}(K_{0}^{(i)},K_{1}^{(i)}){\rm lk}(K_{0}^{(i)},K_{2}^{(i)})
&=&{\rm lk}({K'}_{0}^{(i)},{K'}_{1}^{(i)}){\rm lk}({K'}_{0}^{(i)},{K'}_{2}^{(i)})
\end{eqnarray}
for $i=1,2,\ldots,k$. Then by (\ref{beta_cal2}), (\ref{skeintree2}), (\ref{key0}), (\ref{key1}) and (\ref{key2}), we have that 
\begin{eqnarray}
a_{3}(L')&=&a_{3}(T'_{m})+\sum_{i=1}^{k}\varepsilon(c'_{i})a_{2}(M'_{i}) \nonumber\\
&=&
\sum_{i=1}^{k}\varepsilon(c'_{i})
\Big\{-m{\rm lk}({K'}_{1}^{(i)},{K'}_{2}^{(i)})
+{\rm lk}({K'}_{0}^{(i)},{K'}_{2}^{(i)}){\rm lk}({K'}_{0}^{(i)},{K'}_{1}^{(i)})
\Big\} \nonumber \\
&=&
\sum_{i=1}^{k}\varepsilon(c_{i})
\Big\{-m{\rm lk}(K_{1}^{(i)},K_{2}^{(i)})
+{\rm lk}(K_{0}^{(i)},K_{2}^{(i)}){\rm lk}(K_{0}^{(i)},K_{1}^{(i)})
\Big\}. \label{a3_2}
\end{eqnarray}

On the other hand, by (\ref{skeintree1}) we can see easily that 
\begin{eqnarray}\label{a2a2}
a_{2}(J_{1})+a_{2}(J_{2})
=\sum_{i=1}^{k}\varepsilon(c_{i}){\rm lk}(K_{1}^{(i)},K_{2}^{(i)}). 
\end{eqnarray}
By combining (\ref{a3_1}), (\ref{a3_2}) and (\ref{a2a2}), we have that 
\begin{eqnarray}\label{a3}
a_{3}(L)-a_{3}(L')
&=&\frac{1}{6}(m^{3}-m)
+2m\sum_{i=1}^{k}\varepsilon(c_{i}){\rm lk}(K_{1}^{(i)},K_{2}^{(i)})\nonumber \\
&=&\frac{1}{6}(m^{3}-m)+2m\left\{a_{2}(J_{1})+a_{2}(J_{2})\right\}. 
\end{eqnarray}
Then by (\ref{a3}), we have that 
\begin{eqnarray*}
\tilde{\beta}(L)-\tilde{\beta}(L')
&=&a_{3}(L)-m\left\{a_{2}(J_{1})+a_{2}(J_{2})\right\}
-a_{3}(L')-m\left\{a_{2}(-J_{1})+a_{2}(J_{2})\right\}\\
&=&
a_{3}(L)-a_{3}(L')-2m\left\{a_{2}(J_{1})+a_{2}(J_{2})\right\}\\
&=&
\frac{1}{6}(m^{3}-m). 
\end{eqnarray*}
This completes the proof. 
\end{proof}

\begin{Remark}
{\rm 
Let $f$ be a spatial embedding of a graph $G$ and $\gamma,\ \gamma'$ two disjoint cycles of $G$. By Theorem \ref{gsl_ori}, if ${\rm lk}(f(\gamma),f(\gamma'))=0,\ \pm 1$ then the value of $\tilde{\beta}(f(\gamma),f(\gamma'))$ does not depend on the orientation of $f(\gamma)$ and $f(\gamma')$, namely it is well-defined. But if ${\rm lk}(f(\gamma),f(\gamma'))\neq 0$, then Theorem \ref{gsl_ori} implies that the value of $\tilde{\beta}(f(\gamma),f(\gamma'))$ have the indeterminacy arisen from a choice of the orientations of $f(\gamma)$ and $f(\gamma')$. 
}
\end{Remark}

\section{Definitions of invariants} 

From now onward, we assume that a graph $G$ is {\it oriented}, namely an orientation is given for each edge of $G$. For a subgraph $H$ of $G$, we denote the set of all cycles of $H$ by $\Gamma(H)$. For an edge $e$ of $H$, we denote the set of all oriented cycles of $H$ which contain the edge $e$ and have the orientation induced by the orientation of $e$ by $\Gamma_{e}(H)$. For a pair of two adjacent edges $e$ and $e'$ of $H$, we denote the set of all oriented cycles of $H$ which contain the edges $e$ and $e'$ and have the orientation induced by the orientation of $e$ by $\Gamma_{e,e'}(H)$. We set ${\mathbb Z}_{n}=\{0,1,\ldots,n-1\}$ for a positive integer $n$ and ${\mathbb Z}_{0}={\mathbb Z}$. We call a map $\omega:\Gamma(H)\to {\mathbb Z}_{n}$ a {\it weight on $\Gamma(H)$ over ${\mathbb Z}_{n}$}. Then we say that a weight $\omega$ on $\Gamma(H)$ over ${\mathbb Z}_{n}$ is {\it weakly balanced on an edge $e$} if 
\begin{eqnarray*}
\sum_{\gamma\in\Gamma_{e}(H)}\omega(\gamma)\equiv 0 \pmod{n}, 
\end{eqnarray*}
and {\it weakly balanced on a pair of 
adjacent edges $e$ and $e'$} if 
\begin{eqnarray*}
\sum_{\gamma\in\Gamma_{e,e'}(H)}\omega(\gamma) \equiv 0 \pmod{n}. 
\end{eqnarray*}

Let $G=G_{1}\cup G_{2}$ be a disjoint union of two graphs, $\omega_{i}$ a weight on $\Gamma(G_{i})$ over ${\mathbb Z}_{n}$ $(i=1,2)$ and $f$ a spatial embedding of $G$. Then we say that a weight $\omega_{i}$ is {\it null-homologous on an edge $e$ of $G_{i}$ with respect to $f$ and $\omega_{j}$} ($i\neq j$) if 
\begin{eqnarray*}
{\rm lk}\left(\sum_{\gamma\in \Gamma_{e}(G_{i})}\omega_{i}(\gamma)f(\gamma),f(\gamma')\right) \equiv 0 \pmod{n}
\end{eqnarray*}
for any $\gamma'\in \Gamma(G_{j})$ with $\omega_{j}(\gamma')\neq 0$, and {\it null-homologous on a pair of adjacent edges $e$ and $e'$ of $G_{i}$ with respect to $f$ and $\omega_{j}$} ($i\neq j$) if 
\begin{eqnarray*}
{\rm lk}\left(\sum_{\gamma\in \Gamma_{e,e'}(G_{i})}\omega_{i}(\gamma)f(\gamma),f(\gamma')\right) \equiv 0 \pmod{n}
\end{eqnarray*}
for any $\gamma'\in \Gamma(G_{j})$ with $\omega_{j}(\gamma')\neq 0$. 

\begin{Example}\label{null_homo_ex1}
{\rm 
Let $G=G_{1}\cup G_{2}$ be the graph as illustrated in Fig. \ref{null_homo1}. We denote the cycle $e_{i}\cup e_{j}$ of $G_{1}$ by $\gamma_{ij}$. Let $\omega_{1}$ be the weight on $\Gamma(G_{1})$ over ${\mathbb Z}$ defined by
\begin{eqnarray*}
\omega_{1}(\gamma)=\left\{
   \begin{array}{@{\,}ll}
   1 & (\gamma=\gamma_{12},\ \gamma_{34})\\
   -1 & (\gamma=\gamma_{23},\ \gamma_{14})\\
   0 & ({\rm otherwise}), \\
   \end{array}
\right.
\end{eqnarray*}
and $\omega_{2}$ the weight on $\Gamma(G_{2})$ over ${\mathbb Z}$ defined by $\omega_{2}(\gamma')=1$. Let $f$ be the spatial embedding of $G$ as illustrated in Fig. \ref{null_homo1}. Note that 
\begin{eqnarray*}
\Gamma_{e_{1}}(G_{1})
=\left\{\gamma_{12},\gamma_{13},\gamma_{14}\right\}
=\left\{e_{1}+e_{2},\ e_{1}-e_{3},\ e_{1}+e_{4}\right\}
\end{eqnarray*}
and
\begin{eqnarray*}
\sum_{\gamma\in \Gamma_{e_{1}}(G_{1})}\omega_{1}(\gamma)\gamma
=(e_{1}+e_{2})-(e_{1}+e_{4})=e_{2}-e_{4}.
\end{eqnarray*}
Then we have that 
\begin{eqnarray*}
{\rm lk}\left(\sum_{\gamma\in \Gamma_{e_{1}}(G_{1})}\omega_{1}(\gamma)f(\gamma),f(\gamma')\right) 
={\rm lk}\left(f(e_{2}-e_{4}),f(\gamma')\right)=0.
\end{eqnarray*}
Therefore $\omega_{1}$ is null-homologous on $e_{1}$ with respect to $f$ and $\omega_{2}$. 
\begin{figure}[htbp]
      \begin{center}
\scalebox{0.45}{\includegraphics*{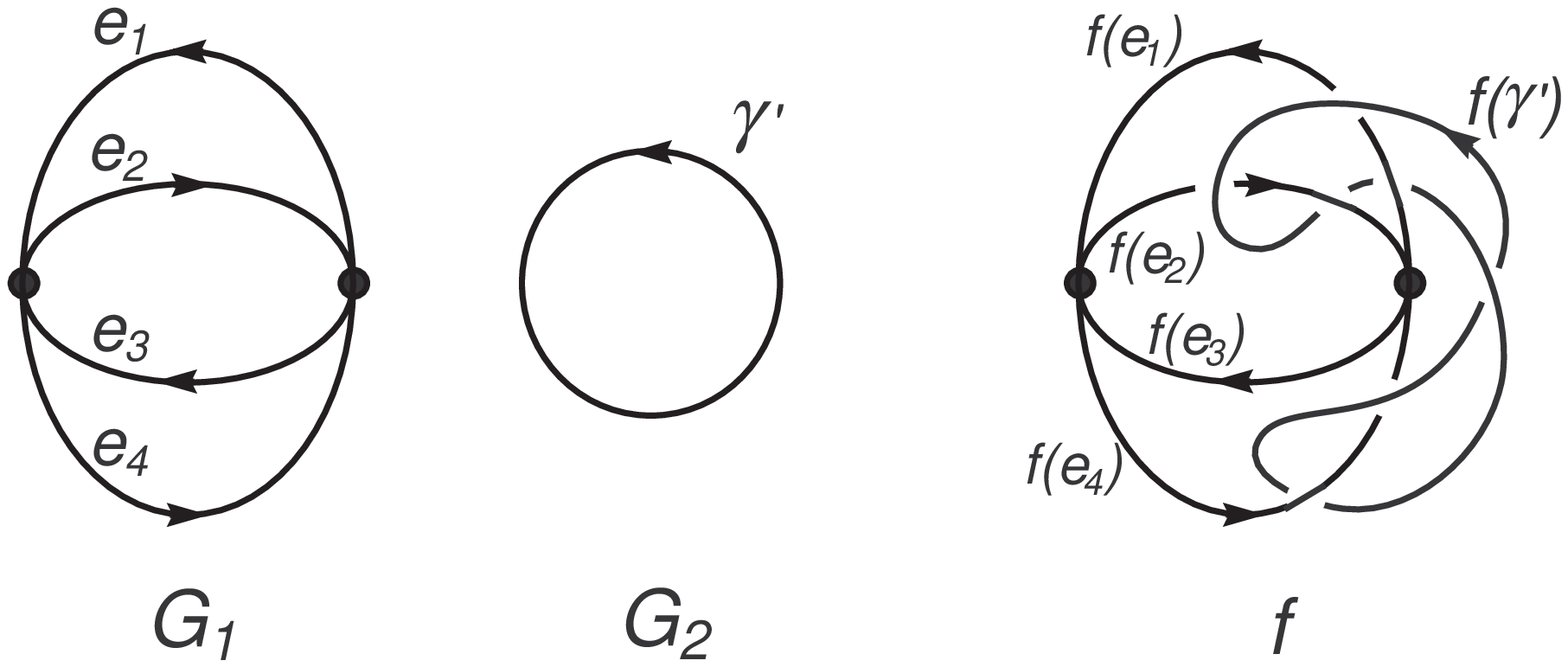}}
      \end{center}
   \caption{}
  \label{null_homo1}
\end{figure} 
}
\end{Example}
\begin{Example}\label{null_homo_ex2}
{\rm 
Let $G=G_{1}\cup G_{2}$ be the graph as illustrated in Fig. \ref{null_homo2}. We denote the cycle of $G_{1}$ which contains $e_{i}$ and $e_{j}$ $(1\le i,j\le 4)$ by $\gamma_{ij}$. Let $\omega_{1}$ be the weight on $\Gamma(G_{1})$ over ${\mathbb Z}$ defined by
\begin{eqnarray*}
\omega_{1}(\gamma)=\left\{
   \begin{array}{@{\,}ll}
   1 & (\gamma=\gamma_{14},\ \gamma_{23})\\
   -1 & (\gamma=\gamma_{13},\ \gamma_{24})\\
   0 & ({\rm otherwise}), \\
   \end{array}
\right.
\end{eqnarray*}
and $\omega_{2}$ the weight on $\Gamma(G_{2})$ over ${\mathbb Z}$ defined by $\omega_{2}(\gamma')=1$. Let $f$ be the spatial embedding of $G$ as illustrated in Fig. \ref{null_homo2}. Note that 
\begin{eqnarray*}
\Gamma_{e_{1},e_{5}}(G_{1})
=\left\{\gamma_{13},\gamma_{14}\right\}
=\left\{e_{1}+e_{5}-e_{3}-e_{6},\ e_{1}+e_{5}-e_{4}-e_{6}\right\}
\end{eqnarray*}
and
\begin{eqnarray*}
\sum_{\gamma\in \Gamma_{e_{1},e_{5}}(G_{1})}\omega_{1}(\gamma)\gamma
=-(e_{1}+e_{5}-e_{3}-e_{6})+(e_{1}+e_{5}-e_{4}-e_{6})=e_{3}-e_{4}.
\end{eqnarray*}
Then we have that 
\begin{eqnarray*}
{\rm lk}\left(\sum_{\gamma\in \Gamma_{e_{1},e_{5}}(G_{1})}\omega_{1}(\gamma)f(\gamma),f(\gamma')\right) 
={\rm lk}\left(f(e_{3}-e_{4}),f(\gamma')\right)=0.
\end{eqnarray*}
Therefore $\omega_{1}$ is null-homologous on a pair of adjacent edges $e_{1}$ and $e_{5}$ with respect to $f$ and $\omega_{2}$. 
\begin{figure}[htbp]
      \begin{center}
\scalebox{0.45}{\includegraphics*{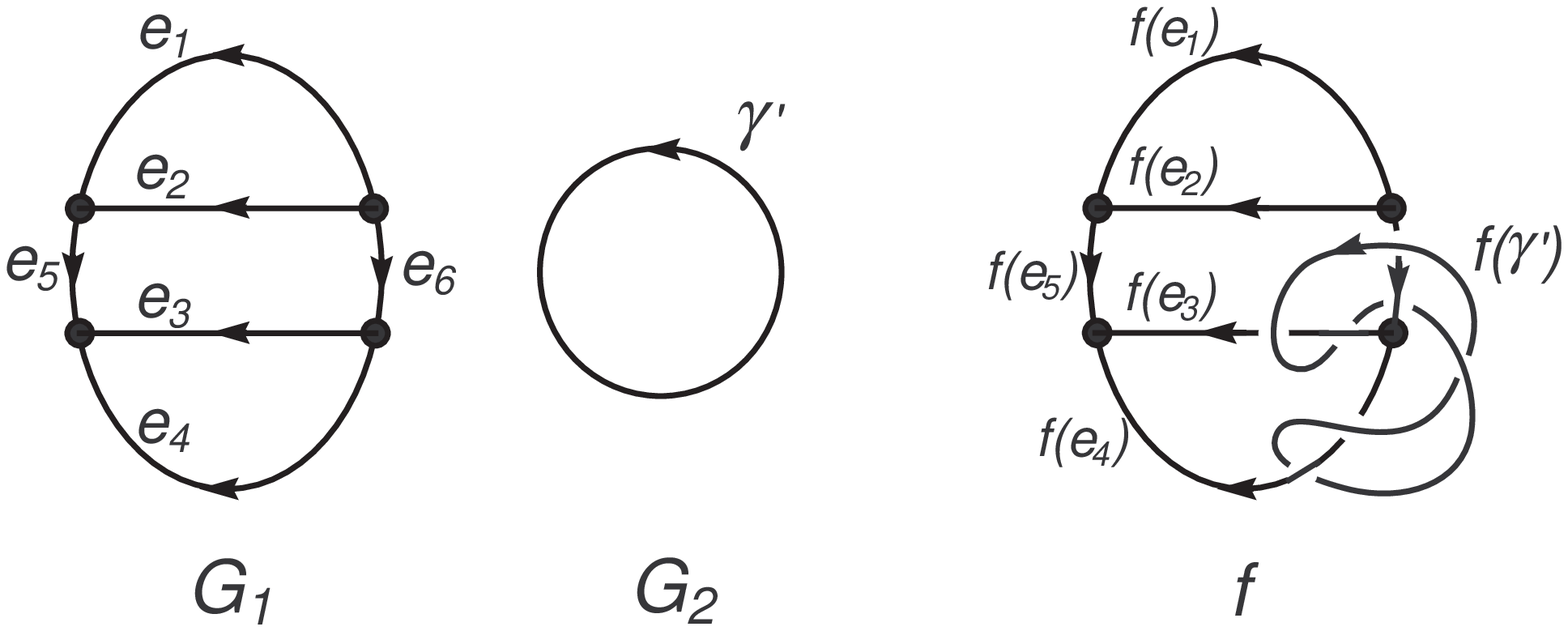}}
      \end{center}
   \caption{}
  \label{null_homo2}
\end{figure} 
}
\end{Example}

Now let $G=G_{1}\cup G_{2}$ be a disjoint union of two graphs, $\omega_{i}$ a weight on $\Gamma(G_{i})$ over ${\mathbb Z}_{n}$ $(i=1,2)$ and $f$ a spatial embedding of $G$. For $\gamma\in \Gamma(G_{1})$ and $\gamma'\in \Gamma(G_{2})$, we put 
\begin{eqnarray*}
\eta(f(\gamma),f(\gamma'))=\frac{1}{6}(m^{3}-m)
\end{eqnarray*}
where $m={\rm lk}(f(\gamma),f(\gamma'))$ under arbitrary orientations of $\gamma$ and $\gamma'$. Then we define that 
\begin{eqnarray*}
\tilde{\eta}_{\omega_{1},\omega_{2}}(f)
={\rm gcd}\left\{
\eta(f(\gamma),f(\gamma'))~|~\gamma\in \Gamma(G_{1}),\ \gamma'\in \Gamma(G_{2}),\ \omega_{1}(\gamma)\omega_{2}(\gamma')\not\equiv 0\pmod{n}
\right\},
\end{eqnarray*}
where ${\rm gcd}$ means the greatest common divisor. Note that $\tilde{\eta}_{\omega_{1},\omega_{2}}(f)$ is a well-defined non-negative integer which does not depends on the choice of orientations of each pair of disjoint cycles. Then we define $\tilde{\beta}_{\omega_{1},\omega_{2}}(f)\in {\mathbb Z}_{n}$ by 
\begin{eqnarray*}
\tilde{\beta}_{\omega_{1},\omega_{2}}(f)\equiv \sum_{\gamma\in\Gamma(G_{1}) \atop \gamma'\in\Gamma(G_{2})}
\omega_{1}(\gamma)\omega_{2}(\gamma')
\tilde{\beta}(f(\gamma),f(\gamma'))\pmod{{\rm gcd}\left\{n,\tilde{\eta}_{\omega_{1},\omega_{2}}(f)\right\}}. 
\end{eqnarray*}
Here we may calculate $\tilde{\beta}(f(\gamma),f(\gamma'))$ under arbitrary orientations of $\gamma$ and $\gamma'$. 

\begin{Remark}
{\rm 
{\rm (1)} For an oriented $2$-component link $L$, $\tilde{\beta}(L)$ is not a link-homotopy invariant of $L$. Thus $\tilde{\beta}_{\omega_{1},\omega_{2}}(f)$ may be not an edge (resp. vertex)-homotopy invariant of $f$ as it is. See also Remark \ref{rem}. 

\noindent
{\rm (2)} By Theorem \ref{gsl_ori}, the value of $\tilde{\beta}(f(\gamma),f(\gamma'))$ is well-defined modulo $\eta(f(\gamma),f(\gamma'))$. This is the reason why we consider the modulo $\tilde{\eta}_{\omega_{1},\omega_{2}}(f)$ reduction. 
}
\end{Remark}

Then, let us state the invariance of $\tilde{\beta}_{\omega_{1},\omega_{2}}$ up to edge (resp. vertex)-homotopy under some conditions on graphs and its spatial embeddings. 

\begin{Theorem}\label{inv1}
{\rm (1)} If $\omega_{i}$ is weakly balanced on any edge of $G_{i}$ and null-homologous on any edge of $G_{i}$ with respect to $f$ and $\omega_{j}$ $(i=1,2,\ i\neq j)$, then $\tilde{\beta}_{\omega_{1},\omega_{2}}(f)$ is an edge-homotopy invariant of $f$. 

\noindent
{\rm (2)} If $\omega_{i}$ is weakly balanced on any pair of adjacent edges 
of $G_{i}$ and null-homologous on any pair of adjacent edges of $G_{i}$ with respect to $f$ and $\omega_{j}$ $(i=1,2,\ i\neq j)$, then $\tilde{\beta}_{\omega_{1},\omega_{2}}(f)$ is a vertex-homotopy invariant of $f$. 
\end{Theorem}

\begin{figure}[htbp]
      \begin{center}
\scalebox{0.4}{\includegraphics*{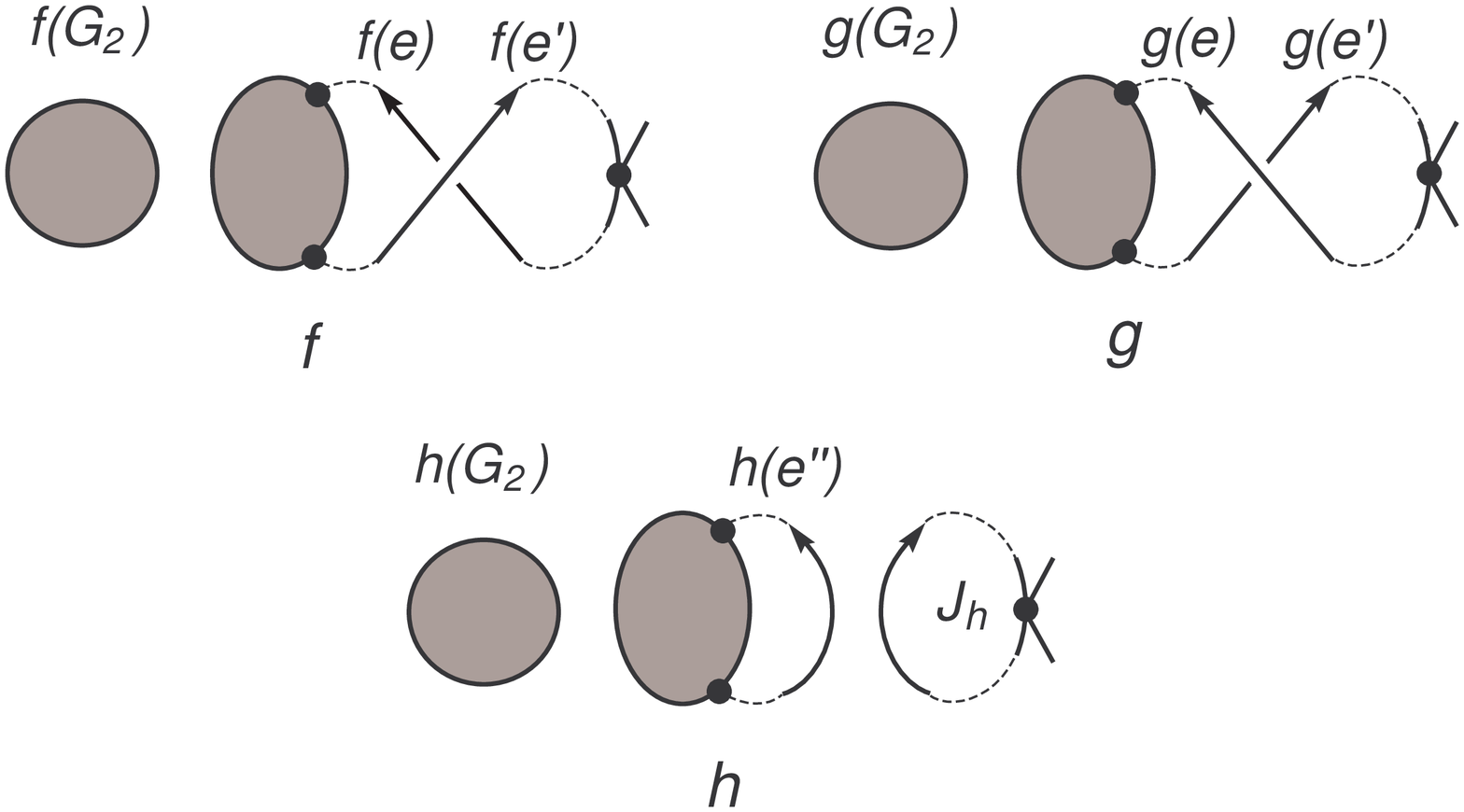}}
      \end{center}
   \caption{}
  \label{skein_sc3}
\end{figure}

\begin{proof}
(1) Let $f$ and $g$ be two spatial embeddings of $G$ such that $g$ is edge-homotopic to $f$. Then it holds that 
\begin{eqnarray}\label{lk0g}
\tilde{\eta}_{\omega_{1},\omega_{2}}(f)=\tilde{\eta}_{\omega_{1},\omega_{2}}(g)
\end{eqnarray}
because the linking number of a constituent $2$-component link of a spatial graph is an edge-homotopy invariant. First we show that if $f$ is transformed into $g$ by self crossing changes on $f(G_{1})$ and ambient isotopies then $\tilde{\beta}_{\omega_{1},\omega_{2}}(f)=\tilde{\beta}_{\omega_{1},\omega_{2}}(g)$. 
It is clear that any link invariant of a constituent link of a spatial graph 
is also an ambient isotopy invariant of the spatial graph. Thus 
we may assume that $g$ is obtained from $f$ by a single crossing change 
on $f(e)$ for an edge $e$ of $G_{1}$ 
as illustrated in Fig. \ref{skein_sc2}. Moreover, by smoothing on this crossing point we can obtain the spatial embedding $h$ of $G$ and the knot $J_{h}$ as illustrated in Fig. \ref{skein_sc2}. Then by (\ref{lk0g}), Lemma \ref{keylemma} and the assumptions for $\omega_{1}$, we have that 
\begin{eqnarray*}
&&\tilde{\beta}_{\omega_{1},\omega_{2}}(f)-\tilde{\beta}_{\omega_{1},\omega_{2}}(g)\\
&\equiv& \sum_{\gamma\in\Gamma_{e}(G_{1}) \atop \gamma'\in\Gamma(G_{2})}
\omega_{1}(\gamma)\omega_{2}(\gamma')
\left\{\tilde{\beta}(f(\gamma),f(\gamma'))-\tilde{\beta}(g(\gamma),g(\gamma'))\right\}\\
&=&\sum_{\gamma\in\Gamma_{e}(G_{1}) \atop \gamma'\in\Gamma(G_{2})}
\omega_{1}(\gamma)\omega_{2}(\gamma'){\rm lk}(h(\gamma'),J_{h})
\left\{{\rm lk}(f(\gamma),f(\gamma'))-{\rm lk}(h(\gamma'),J_{h})\right\}\\
&=& \sum_{\gamma'\in \Gamma(G_{2})}\omega_{2}(\gamma')
\Bigg\{
{\rm lk}(h(\gamma'),J_{h})\sum_{\gamma\in \Gamma_{e}(G_{1})}\omega_{1}(\gamma){\rm lk}(f(\gamma),f(\gamma'))\\
&&-\sum_{\gamma\in \Gamma_{e}(G_{1})}\omega_{1}(\gamma){\rm lk}(h(\gamma'),J_{h})^{2}
\Bigg\}\\
&=&
\sum_{\gamma'\in \Gamma(G_{2})}\omega_{2}(\gamma')
\Bigg\{
{\rm lk}(h(\gamma'),J_{h}){\rm lk}\left(\sum_{\gamma\in \Gamma_{e}(G_{1})}\omega_{1}(\gamma)f(\gamma),f(\gamma')\right)\\
&&-{\rm lk}(h(\gamma'),J_{h})^{2}\left(\sum_{\gamma\in \Gamma_{e}(G_{1})}\omega_{1}(\gamma)\right)
\Bigg\}\\
&\equiv&0\pmod{{\rm gcd}\left\{n,\tilde{\eta}_{\omega_{1},\omega_{2}}(f)\right\}}. 
\end{eqnarray*}
Therefore we have that $\tilde{\beta}_{\omega_{1},\omega_{2}}(f)=\tilde{\beta}_{\omega_{1},\omega_{2}}(g)$.  In the same way we can show that if $f$ is transformed into $g$ by 
self crossing changes on $f(G_{2})$ and ambient isotopies then 
$\tilde{\beta}_{\omega_{1},\omega_{2}}(f)= \tilde{\beta}_{\omega_{1},\omega_{2}}(g)$. Thus we have that $\tilde{\beta}_{\omega_{1},\omega_{2}}(f)$ is an 
edge-homotopy invariant of $f$. 
\begin{figure}[htbp]
      \begin{center}
\scalebox{0.4}{\includegraphics*{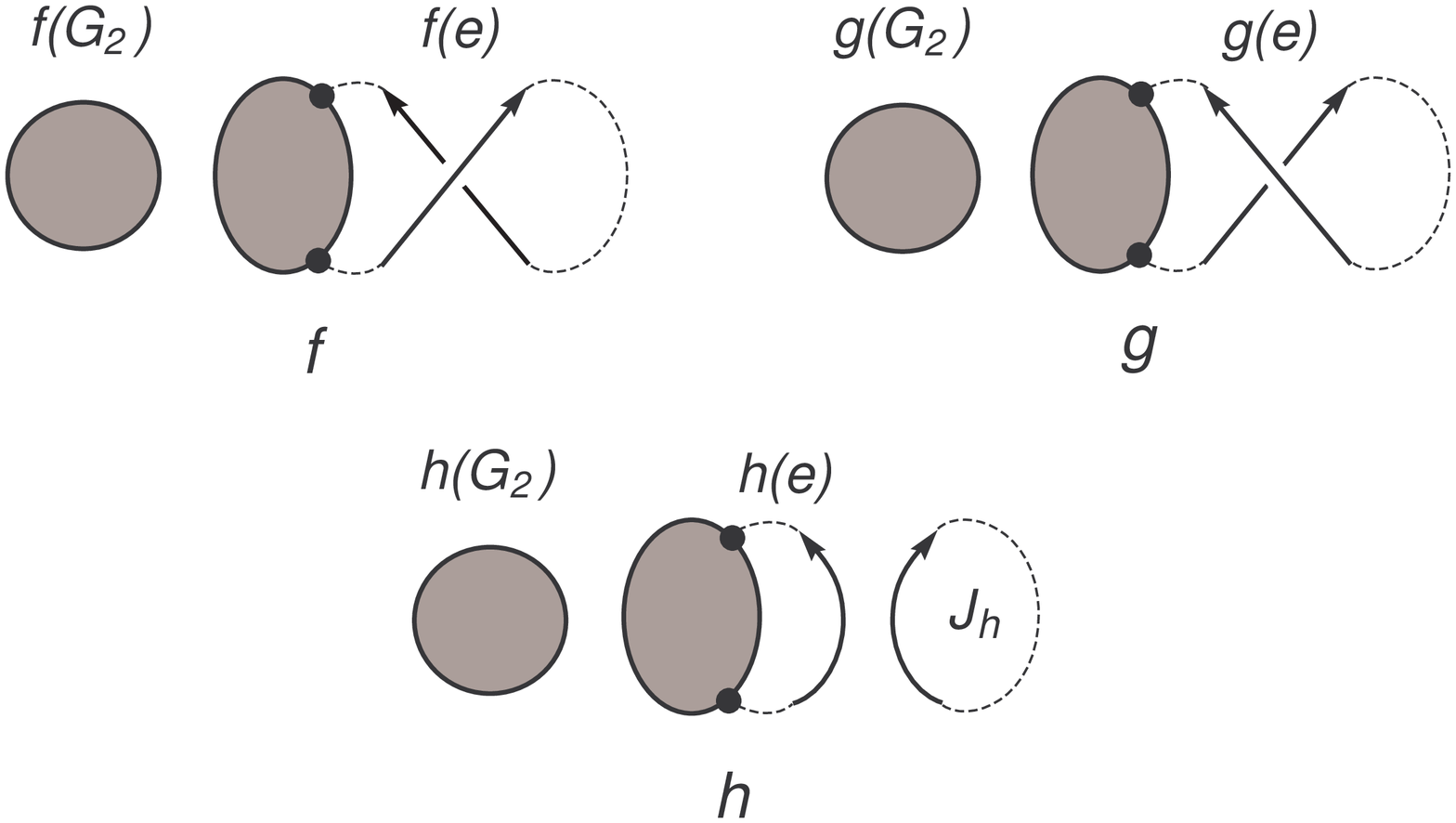}}
      \end{center}
   \caption{}
  \label{skein_sc2}
\end{figure} 

(2) By considering 
the triple of spatial embeddings as illustrated in Fig. \ref{skein_sc3}, we can prove (2) in a similar way as the proof of (1). 
We omit the details. 
\end{proof}
\begin{Remark}\label{gene_fn}
{\rm 
In particular, if it holds that 
\begin{eqnarray*}
\omega_{1}(\gamma)\omega_{2}(\gamma'){\rm lk}(f(\gamma),f(\gamma'))=0
\end{eqnarray*}
for any $\gamma\in \Gamma(G_{1})$ and $\gamma'\in \Gamma(G_{2})$, then $\tilde{\beta}_{\omega_{1},\omega_{2}}(f)$ coincides with Fleming and the author's invariant ${\beta}_{\omega_{1},\omega_{2}}(f)$ defined in \cite{fleming-nikkuni05}. 
}
\end{Remark}
\section{Examples} 

Let $G$ be a planar graph which is not a cycle. An embedding $p:G\to {\mathbb S}^{2}$ is said to be {\it cellular} if the closure of each of the connected components of ${\mathbb S^{2}-p(G)}$ on ${\mathbb S}^{2}$ is homeomorphic to the disk. Then we regard the set of the boundaries of all of the connected components of ${\mathbb S}^{2}-p(G)$ as a subset of $\Gamma(G)$ and denote it by $\Gamma_{p}(G)$. We say that {\it $G$ admits a checkerboard coloring on ${\mathbb S}^{2}$} if there exists a cellular embedding $p:G\to {\mathbb S}^{2}$ such that we can color all of the connected components of ${\mathbb S}^{2}-p(G)$ by two colors (black and white) so that any of the two components which are adjacent by at least one edge have distinct colors. We denote the subset of $\Gamma_{p}(G)$ which corresponds to the black (resp. white) colored components by $\Gamma_{p}^{b}(G)$ (resp. $\Gamma_{p}^{w}(G)$). Then, for any edge $e$ of $G$, there exist exactly two cycles $\gamma\in \Gamma_{p}^{b}(G)$ and $\gamma'\in \Gamma_{p}^{w}(G)$ such that $e\subset\gamma$ and $e\subset\gamma'$. Thus we have the following immediately. 

\begin{Proposition}\label{checkerboard_weight}
Let $G$ be a planar graph which is not a cycle and admits a checkerboard coloring on ${\mathbb S}^{2}$ with respect to a cellular embedding $p:G\to {\mathbb S}^{2}$. Let $\omega_{p}$ be the weight on $\Gamma(G)$ over ${\mathbb Z}_{n}$ defined by 
\begin{eqnarray*}
\omega_{p}(\gamma)= \left\{
\begin{array}{@{\,}ll}
1 & \mbox{$(\gamma\in \Gamma_{p}^{b}(G))$} \\
n-1 &   \mbox{$(\gamma\in\Gamma_{p}^{w}(G))$} \\
0 &   \mbox{$(\gamma\in \Gamma(G)-\Gamma_{p}(G))$.}
\end{array}
\right. 
\end{eqnarray*}
Then $\omega_{p}$ is weakly balanced on any edge of $G$. 
\end{Proposition}
We call the weight $\omega_{p}$ in Proposition \ref{checkerboard_weight} a 
{\it checkerboard weight}. Moreover, by giving the counter clockwise orientation to each $p(\gamma)$ for $\gamma\in \Gamma_{p}^{b}(G)$ and the clockwise orientation to each $p(\gamma)$ for $\gamma\in \Gamma_{p}^{w}(G)$ with respect to the orientation of ${\mathbb S}^{2}$, an orientation is given for each edge of $G$ naturally. We call this orientation of $G$ a {\it checkerboard orientation} over the checkerboard coloring. Since the orientation of each edge $e$ is coherent with the orientation of each cycle $\gamma\in \Gamma_{p}(G)$ which contains $e$, by Theorem \ref{inv1} we have the following. 

\begin{Theorem}\label{integer_inv}
Let $G=G_{1}\cup G_{2}$ be a disjoint union of two planar graphs such that $G_{i}$ is not a cycle and admits a checkerboard coloring on ${\mathbb S}^{2}$ with respect to a cellular embedding $p_{i}:G\to {\mathbb S}^{2}$ $(i=1,2)$. Let $\omega_{p_{i}}$ be the checkerboard weight on $\Gamma(G_{i})$ over ${\mathbb Z}_{n}$ $(i=1,2)$. We orient $G$ by the checkerboard orientation of $G_{i}$ over the checkerboard coloring $(i=1,2)$. Then, for a spatial embedding $f$ of $G$, if $\omega_{i}$ is null-homologous on any edge of $G_{i}$ with respect to $f$ and $\omega_{j}$ $(i=1,2,\ i\neq j)$, then $\tilde{\beta}_{\omega_{1},\omega_{2}}(f)\pmod{n}$ is an edge-homotopy invariant of $f$. 
\end{Theorem}
\begin{Example}\label{integer_inv_ex}
{\rm 
Let $G=G_{1}\cup G_{2}$ be a disjoint union of two planar graphs as in Theorem \ref{integer_inv} and $f$ a spatial embedding of $G$. Let $\omega_{p_{i}}:\Gamma(G_{i})\to {\mathbb Z}_{n}$ be the checkerboard weight ($i=1,2$), where 
\begin{eqnarray*}
n={\rm gcd}\left\{{\rm lk}(f(\gamma),f(\gamma'))~|~\gamma\in \Gamma_{p_{1}}(G_{1}),\ \gamma'\in \Gamma_{p_{2}}(G_{2})\right\}.
\end{eqnarray*}
Then, for any edge $e$ of $G_{i}$ and any $\gamma'\in \Gamma_{p_{j}}(G_{j})$ $(i\neq j)$, we have that 
\begin{eqnarray*}
{\rm lk}\left(\sum_{\gamma\in \Gamma_{e}(G_{i})}\omega_{i}(\gamma)f(\gamma),f(\gamma')
\right)
=
\sum_{\gamma\in \Gamma_{e}(G_{i})}\omega_{i}(\gamma){\rm lk}\left(f(\gamma),f(\gamma')\right)
\equiv 0\pmod{n}. 
\end{eqnarray*}
Thus we have that $\omega_{i}$ is null-homologous on any edge of $G_{i}$ with respect to $f$ and $\omega_{j}$ $(i=1,2,\ i\neq j)$. Therefore we have that $\tilde{\beta}_{\omega_{p_{1}},\omega_{p_{2}}}(f)\pmod{n}$ is an edge-homotopy invariant of $f$. 

For example, let $\Theta_{4}$ be the graph with two vertices $u$ and $v$ and $4$ edges $e_{1},e_{2},e_{3},e_{4}$ each of which joins $u$ and $v$. We denote the cycle of $\Theta_{4}$ consists of two edges $e_{i}$ and $e_{j}$ by $\gamma_{ij}$. Let $p:\Theta_{4}\to {\mathbb S}^{2}$ be the cellular embedding as illustrated in the left-hand side of Fig. \ref{theta_4_check}. It is clear that $\Theta_{4}$ admits the checkerboard coloring on ${\mathbb S}^{2}$ with respect to $p$ as illustrated in the center of Fig. \ref{theta_4_check}. The right-hand side of Fig. \ref{theta_4_check} shows the checkerboard orientation of $\Theta_{4}$ over the checkerboard coloring. 
\begin{figure}[htbp]
      \begin{center}
\scalebox{0.45}{\includegraphics*{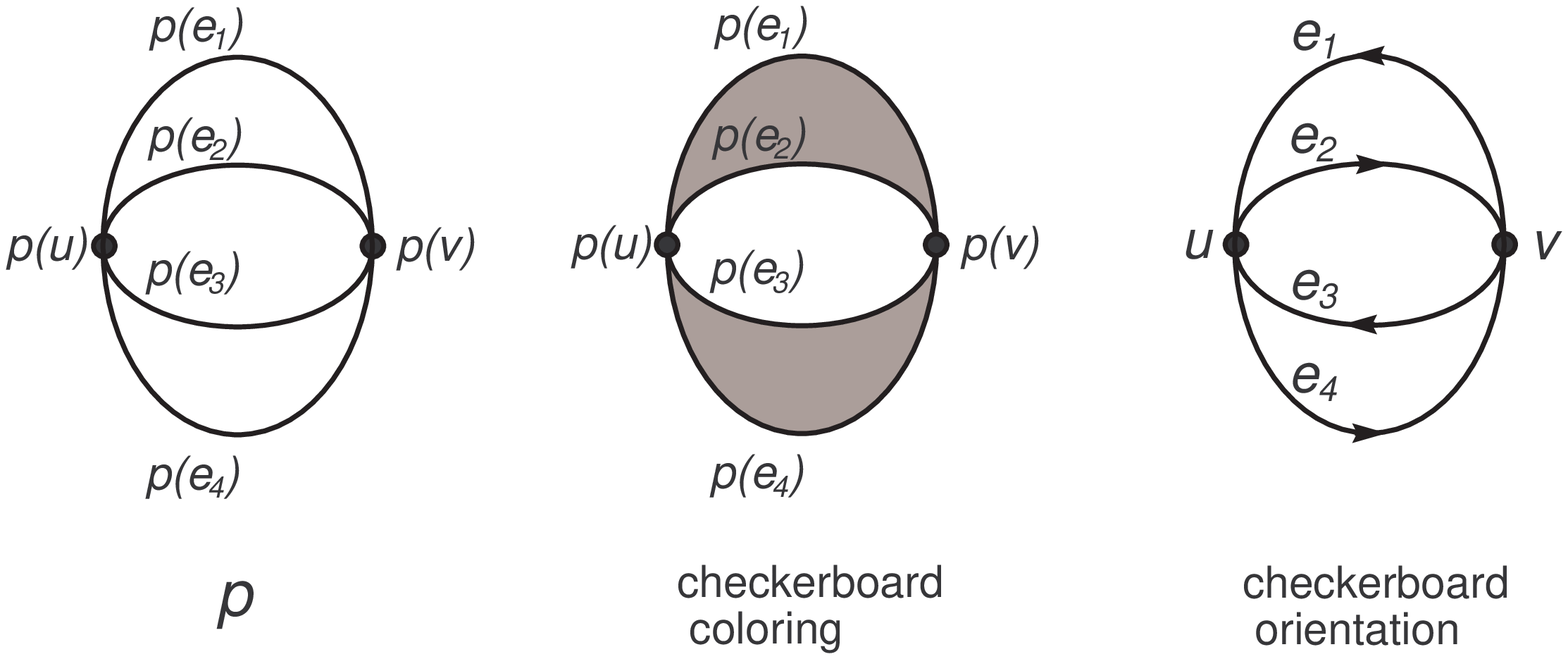}}
      \end{center}
   \caption{}
  \label{theta_4_check}
\end{figure} 

Let $G=\Theta_{4}^{1}\cup \Theta_{4}^{2}$ be a disjoint union of two copies of $\Theta_{4}$. For a non-negative integer $m$, let $f_{m}$ and $g_{m}$ be two spatial embeddings of $G$ as illustrated in Fig. \ref{newex1}. Note that 
\begin{eqnarray*}
{\rm lk}(f_{m}(\gamma),f_{m}(\gamma'))={\rm lk}(g_{m}(\gamma),g_{m}(\gamma'))=0\ {\rm or}\ m
\end{eqnarray*}
 for any $\gamma\in \Gamma(\Theta_{4}^{1})$ and $\gamma' \in \Gamma(\Theta_{4}^{2})$. So we have that $n=m$. Let $\omega_{i}:\Gamma(\Theta_{4}^{i})\to {\mathbb Z}_{m}$ be the checkerboard weight ($i=1,2$). Then, by a direct calculation we can see that the constituent $2$-component link of $f_{m}$ which has a non-zero generalized Sato-Levine invariant is only $L=f_{m}(\gamma_{14}\cup \gamma'_{14})$ and $\tilde{\beta}(L)=2$. Actually the other constituent $2$-component link $f_{m}(\gamma\cup \gamma')$ for $\gamma\in \Gamma_{p}(\Theta_{4}^{1})$ and $\gamma' \in \Gamma_{p}(\Theta_{4}^{2})$ is a trivial $2$-component link or $T'_{m}$ as illustrated in Fig. \ref{(2,m)torus}. Thus we have that $\tilde{\beta}_{\omega_{1},\omega_{2}}(f_{m})\equiv 2\pmod{m}$. On the other hand, we can see that each constituent $2$-component link $g_{m}(\gamma\cup \gamma')$ for $\gamma\in \Gamma_{p}(\Theta_{4}^{1})$ and $\gamma' \in \Gamma_{p}(\Theta_{4}^{2})$ is a trivial $2$-component link or $T'_{m}$. Thus we have that $\tilde{\beta}_{\omega_{1},\omega_{2}}(g_{m})\equiv 0\pmod{m}$. Therefore we have that $f_{m}$ and $g_{m}$ are not edge-homotopic if $m\neq 1,2$. We remark here that the case of $m=0$ has already shown by Fleming and the author in \cite[Example 4.3]{fleming-nikkuni05}. 
\begin{figure}[htbp]
      \begin{center}
\scalebox{0.335}{\includegraphics*{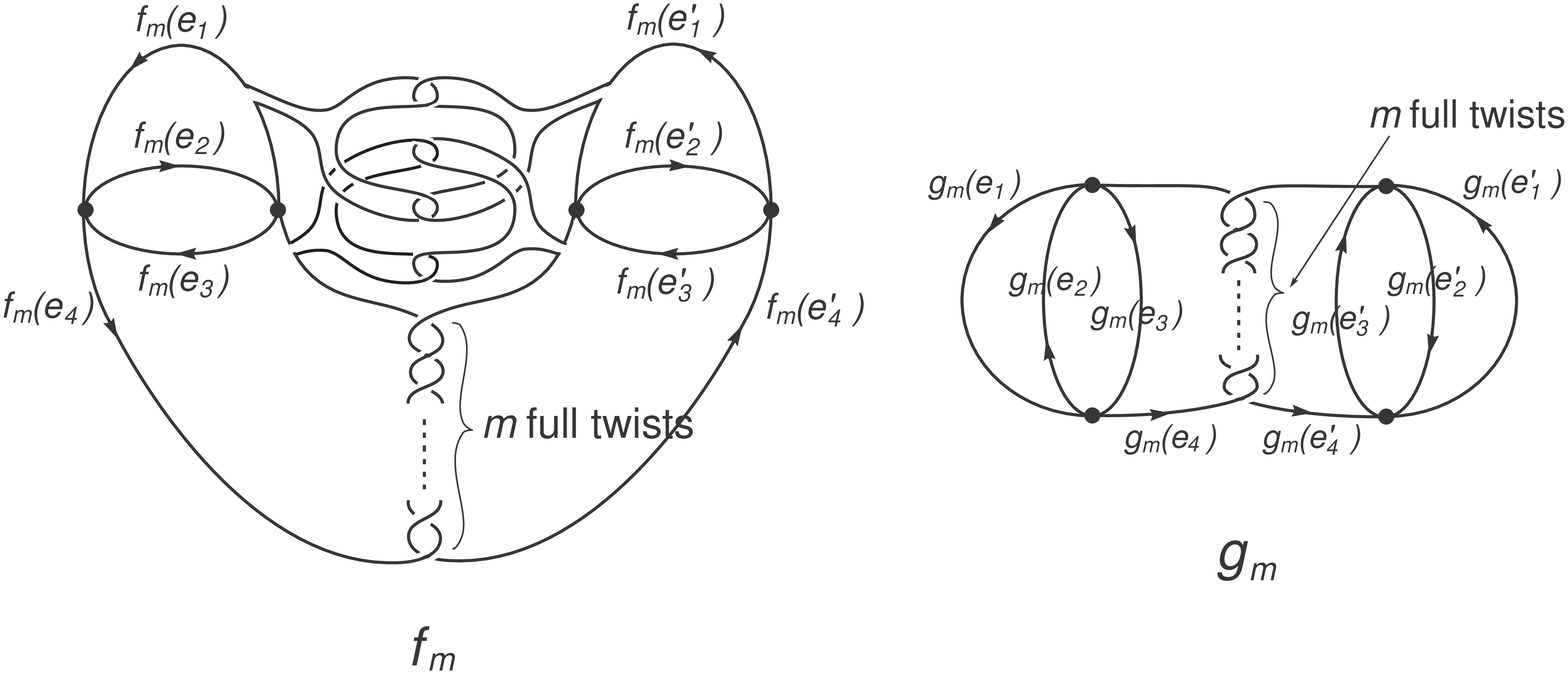}}
      \end{center}
   \caption{}
  \label{newex1}
\end{figure} 
}
\end{Example}
\begin{Example}\label{integer_inv_ex2}
{\rm 
Let $G=\Theta_{4}^{1}\cup \Theta_{4}^{2}$ be a disjoint union of two copies of $\Theta_{4}$ oriented in the same way as Example \ref{integer_inv_ex} and $\omega_{i}:\Gamma(\Theta_{4}^{i})\to {\mathbb Z}$ the checkerboard weight ($i=1,2$). For $\Theta_{4}^{1}$, we have that 
\begin{eqnarray*}
&&\sum_{\gamma\in \Gamma_{e_{1}}(\Theta_{4}^{1})}\omega_{1}(\gamma)\gamma=e_{2}-e_{4},\ \ 
\sum_{\gamma\in \Gamma_{e_{2}}(\Theta_{4}^{1})}\omega_{1}(\gamma)\gamma=e_{1}-e_{3},\\
&&\sum_{\gamma\in \Gamma_{e_{3}}(\Theta_{4}^{1})}\omega_{1}(\gamma)\gamma=e_{4}-e_{2},\ \ 
\sum_{\gamma\in \Gamma_{e_{4}}(\Theta_{4}^{1})}\omega_{1}(\gamma)\gamma=e_{3}-e_{1}.\
\end{eqnarray*}
This implies that $\omega_{1}$ is null-homologous on any edge of $G_{1}$ with respect to a spatial embedding $f$ of $G$ and $\omega_{2}$ if and only if 
\begin{eqnarray}\label{condition1}
{\rm lk}\left(
f(\gamma_{13}),f(\gamma'))
\right)=
{\rm lk}\left(
f(\gamma_{24}),f(\gamma'))
\right)=0
\end{eqnarray}
for any $\gamma'\in \Gamma_{p}(\Theta_{4}^{2})$. The same condition can be said of $\omega_{2}$. For an integer $m$, let $f_{m}$ be the spatial embedding of $G$ as illustrated in Fig. \ref{newex2}. Note that 
\begin{eqnarray*}
{\rm lk}(f_{k}(\gamma),f_{k}(\gamma'))={\rm lk}(f_{l}(\gamma),f_{l}(\gamma'))=0\ {\rm or}\ 1\ (k\neq l) 
\end{eqnarray*}
 for any $\gamma\in \Gamma(\Theta_{4}^{1})$ and $\gamma' \in \Gamma(\Theta_{4}^{2})$. Since we can see that $\omega_{i}$ satisfies (\ref{condition1}), we have that $\omega_{i}$ is null-homologous on any edge of $G_{i}$ with respect to $f_{m}$ and $\omega_{j}$ $(i=1,2,\ i\neq j)$. Namely $\tilde{\beta}_{\omega_{1},\omega_{2}}(f_{m})$ is an integer-valued edge-homotopy invariant of $f_{m}$. Then, by a direct calculation we can see that the constituent $2$-component link of $f_{m}$ which has a non-zero generalized Sato-Levine invariant is only $L=f_{m}(\gamma_{14}\cup \gamma'_{14})$ and $\tilde{\beta}(L)=2m$. Actually the other constituent $2$-component link $f_{m}(\gamma\cup \gamma')$ for $\gamma\in \Gamma_{p}(\Theta_{4}^{1})$ and $\gamma' \in \Gamma_{p}(\Theta_{4}^{2})$ is a Hopf link. Thus we have that $\tilde{\beta}_{\omega_{1},\omega_{2}}(f_{m})=2m$. Therefore we have that $f_{k}$ and $f_{l}$ are not edge-homotopic for $k\neq l$. 
\begin{figure}[htbp]
      \begin{center}
\scalebox{0.4}{\includegraphics*{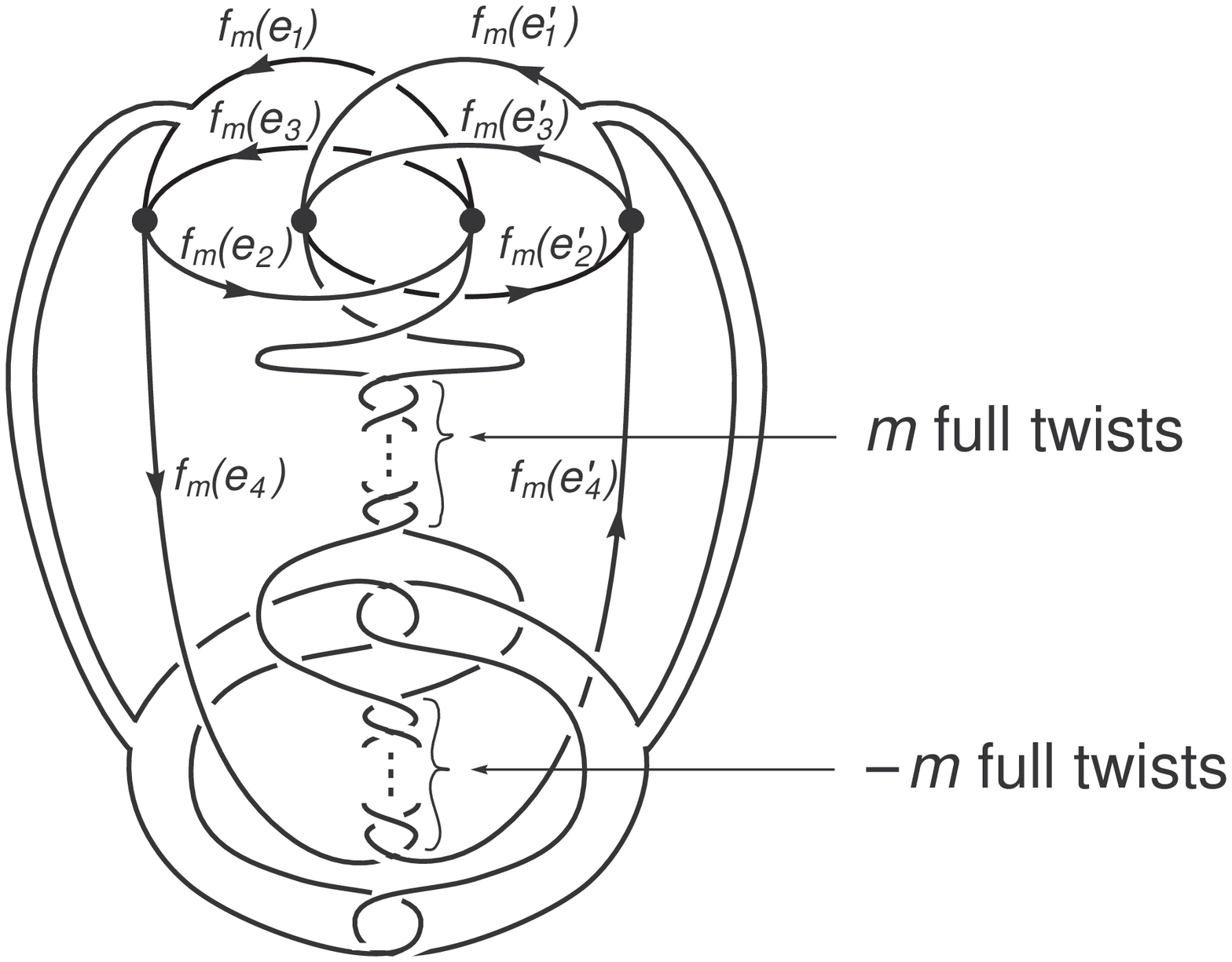}}
      \end{center}
   \caption{}
  \label{newex2}
\end{figure} 
}
\end{Example}
\begin{Example}\label{integer_inv_ex3}
{\rm 
Let $H$ be the same oriented graph as $G_{1}$ in Fig. \ref{null_homo2}. We denote the cycle of $H$ which contains $e_{i}$ and $e_{j}$ by $\gamma_{ij}$. Let $G=H^{1}\cup H^{2}$ be a disjoint union of two copies of $H$ and $f_{m}$ the spatial embedding of $G$ as illustrated in Fig. \ref{newex3}. This spatial embedding $f_{m}$ contains exactly one constituent $4$-component link $L_{m}=f_{m}(\gamma_{12}\cup \gamma_{34}\cup \gamma'_{12}\cup \gamma'_{34})$. By calculating Milnor's {\it $\mu$-invariant} \cite{milnor54} of length $4$ of $L_{m}$, it can be shown that $f_{k}$ and $f_{l}$ are not vertex-homotopic for $k\neq l$. But we can also prove this fact by our vertex-homotopy invariant as follows. Let $\omega$ be the same weight on $\Gamma(H)$ over ${\mathbb Z}$ as $\omega_{1}$ in Example \ref{null_homo_ex2}. We define the weight $\omega_{i}$ of $\Gamma(H^{i})$ over ${\mathbb Z}$ in the same way as $\omega$ ($i=1,2$). It is easy to see that $\omega_{i}$ is weakly balanced on any pair of adjacent edges of $H^{i}$ ($i=1,2$). Moreover, we have that 
\begin{eqnarray*}
&&\sum_{\gamma\in \Gamma_{e_{1},e_{5}}(H)}\omega(\gamma)\gamma=e_{3}-e_{4},\ \ 
\sum_{\gamma\in \Gamma_{e_{1},e_{2}}(H)}\omega(\gamma)\gamma=0,\\
&&\sum_{\gamma\in \Gamma_{e_{1},e_{6}}(H)}\omega(\gamma)\gamma=e_{3}-e_{4},\ \ 
\sum_{\gamma\in \Gamma_{e_{2},e_{5}}(H)}\omega(\gamma)\gamma=e_{4}-e_{3},\\
&&\sum_{\gamma\in \Gamma_{e_{2},e_{6}}(H)}\omega(\gamma)\gamma=e_{4}-e_{3},\ \ 
\sum_{\gamma\in \Gamma_{e_{3},e_{5}}(H)}\omega(\gamma)\gamma=e_{1}-e_{2},\\
&&\sum_{\gamma\in \Gamma_{e_{3},e_{6}}(H)}\omega(\gamma)\gamma=e_{1}-e_{2},\ \ 
\sum_{\gamma\in \Gamma_{e_{3},e_{4}}(H)}\omega(\gamma)\gamma=0,\\
&&\sum_{\gamma\in \Gamma_{e_{4},e_{5}}(H)}\omega(\gamma)\gamma=e_{2}-e_{1},\ \ 
\sum_{\gamma\in \Gamma_{e_{4},e_{6}}(H)}\omega(\gamma)\gamma=e_{2}-e_{1}.
\end{eqnarray*}
This implies that $\omega_{1}$ is null-homologous on any pair of adjacent edges of $H^{1}$ with respect to a spatial embedding $f$ of $G$ and $\omega_{2}$ if and only if 
\begin{eqnarray}\label{condition3}
{\rm lk}\left(
f(\gamma_{12}),f(\gamma'))
\right)=
{\rm lk}\left(
f(\gamma_{34}),f(\gamma'))
\right)=0
\end{eqnarray}
for any $\gamma'\in \Gamma(H^{2})$ and $\omega_{2}(\gamma')\neq 0$. The same condition can be said of $\omega_{2}$. Then we can see that $\omega_{i}$ satisfies (\ref{condition3}) for $f_{m}$, namely $\omega_{i}$ is null-homologous on any edge of $G_{i}$ with respect to $f_{m}$ and $\omega_{j}$ $(i=1,2,\ i\neq j)$. Namely $\tilde{\beta}_{\omega_{1},\omega_{2}}(f_{m})$ is a vertex-homotopy invariant of $f_{m}$. Since the linking number of each constituent $2$-component link of $f_{m}$ is $0$ or $\pm 1$, we have that $\tilde{\eta}_{\omega_{1},\omega_{2}}(f_{m})=0$. Namely $\tilde{\beta}_{\omega_{1},\omega_{2}}(f_{m})$ is integer-valued. By a direct calculation we have that $\tilde{\beta}_{\omega_{1},\omega_{2}}(f_{m})=2m$ in the same way as Example \ref{integer_inv_ex2}. Therefore we have that $f_{k}$ and $f_{l}$ are not vertex-homotopic for $k\neq l$. 
\begin{figure}[htbp]
      \begin{center}
\scalebox{0.4}{\includegraphics*{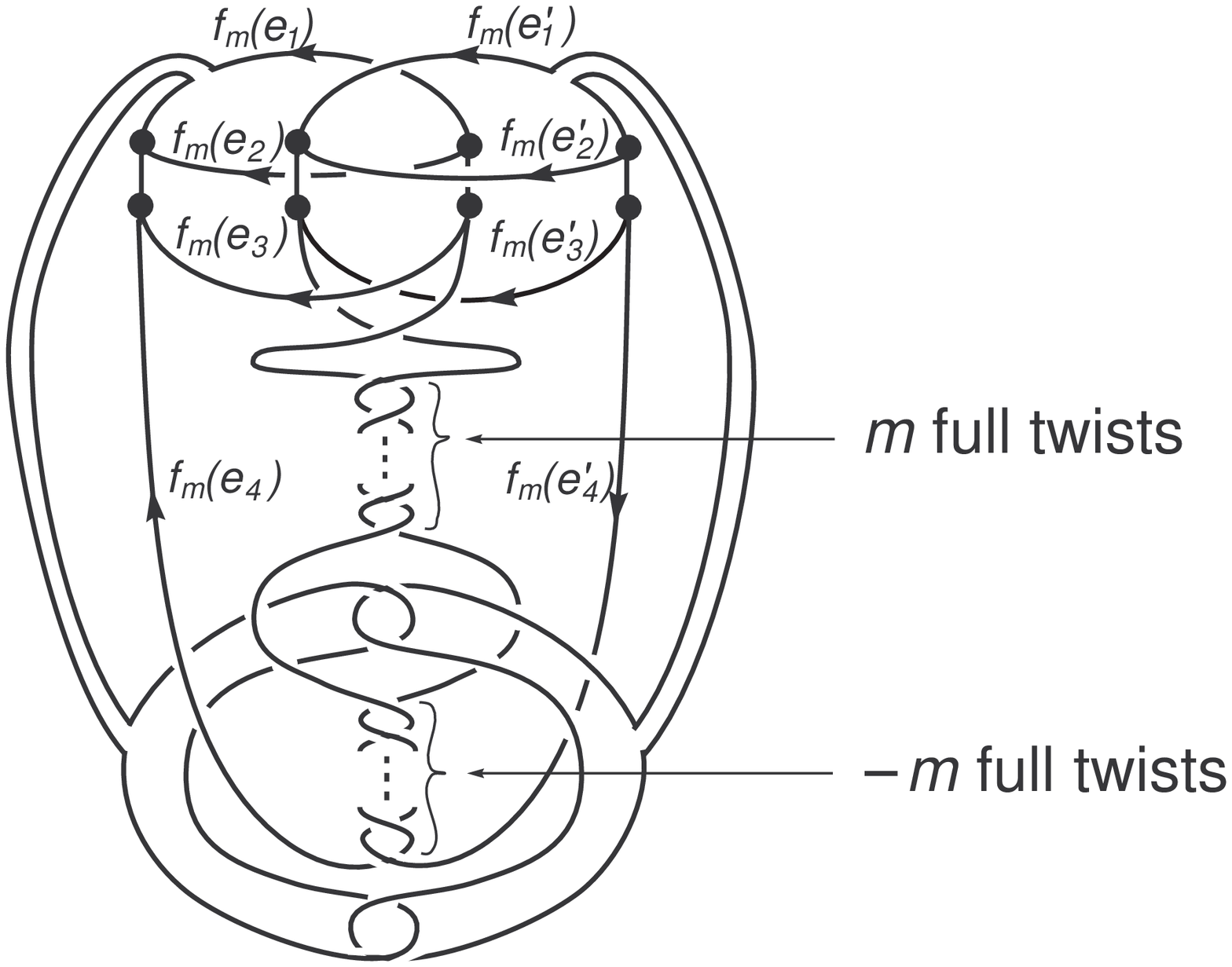}}
      \end{center}
   \caption{}
  \label{newex3}
\end{figure} 
}
\end{Example}

\begin{Remark}\label{rem}
{\rm In Theorems \ref{inv1} and \ref{integer_inv}, the condition ``$\omega_{i}$ is null-homologous on any edge of $G_{i}$ with respect to $f$ and $\omega_{j}$ $(i=1,2,\ i\neq j)$'' is essential. Let $G=\Theta_{4}^{1}\cup \Theta_{4}^{2}$ be a disjoint union of two copies of $\Theta_{4}$ oriented in the same way as Example \ref{integer_inv_ex2} and $\omega_{i}:\Gamma(\Theta_{4}^{i})\to {\mathbb Z}$ the checkerboard weight ($i=1,2$). Let $f$ and $g$ be two spatial embeddings of $G$ as illustrated in Fig. \ref{null_homo_remark}. Note that $f$ and $g$ are edge-homotopic. But by a direct calculation we have that $\tilde{\beta}_{\omega_{1},\omega_{2}}(f)=-1$ and $\tilde{\beta}_{\omega_{1},\omega_{2}}(g)=0$, namely $\tilde{\beta}_{\omega_{1},\omega_{2}}(f)$ is not an edge-homotopy invariant of $f$. Actually $\omega_{1}$ is not null-homologous on $e_{4}$ with respect to $f$ and $\omega_{2}$. 
\begin{figure}[htbp]
      \begin{center}
\scalebox{0.4}{\includegraphics*{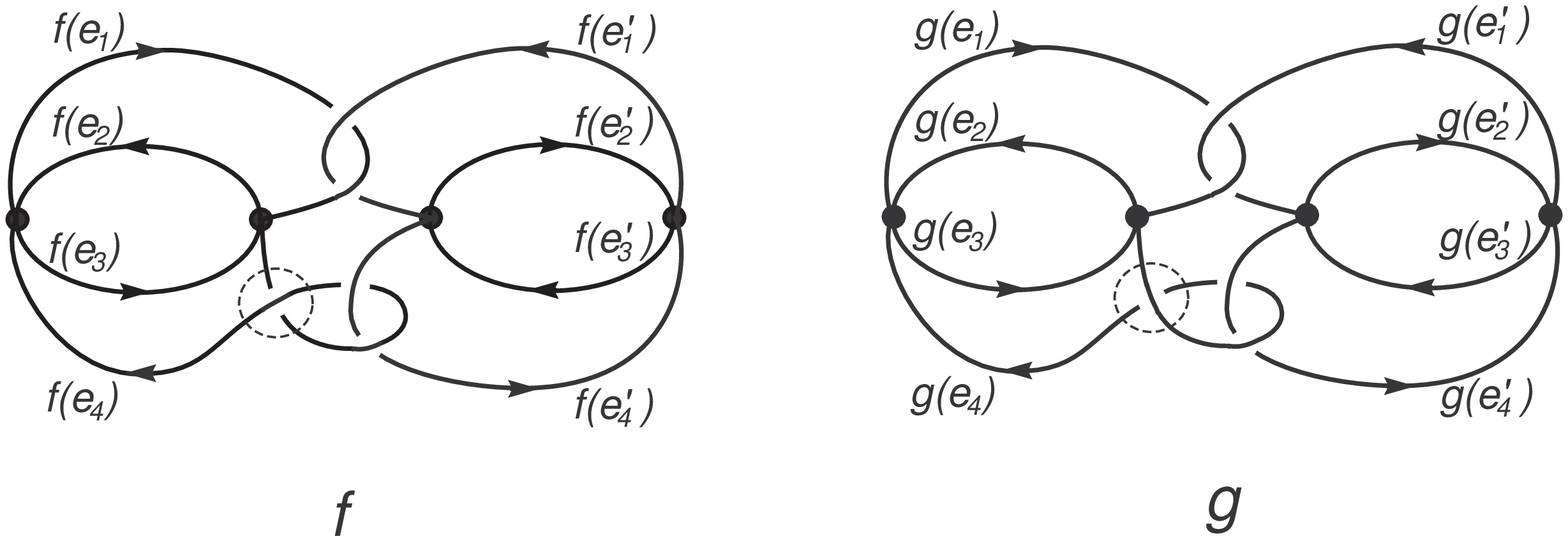}}
      \end{center}
   \caption{}
  \label{null_homo_remark}
\end{figure} 
}
\end{Remark}

\section*{Acknowledgment}

The author would like to thank Professor Kouki Taniyama for informing him about the results in \cite{taniyama0}. He is also grateful to Doctor Thomas Fleming for his valuable comments. 

{\normalsize
}

\end{document}